# UNFOLDING OF CHAOTIC QUADRATIC MAPS — PARAMETER DEPENDENCE OF NATURAL MEASURES

HANS THUNBERG

ABSTRACT. We consider perturbations of quadratic maps $f_a$ admitting an absolutely continuous invariant probability measure, where $a$ is in a certain positive measure set $\mathcal{A}$ of parameters, and show that in any neighborhood of any such an $f_a$, we find a rich fauna of dynamics. There are maps with periodic attractors as well as non-periodic maps whose critical orbit is absorbed by the continuation of any prescribed hyperbolic repeller of $f_a$. In particular, Misiurewicz maps are dense in $\mathcal{A}$.

Almost all maps $f_a$ in the quadratic family is known to possess a unique natural measure, that is, an invariant probability measure $\mu_a$ describing the asymptotic distribution of almost all orbits. We discuss weak*-(dis)continuity properties of the map $a \mapsto \mu_a$ near the set $\mathcal{A}$, and prove that almost all maps in $\mathcal{A}$ have the property that $\mu_a$ can be approximated with measures supported on periodic attractors of certain nearby maps. On the other hand, for any $a \in \mathcal{A}$ and any periodic repeller $\Gamma_a$ of $f_a$, the singular measure supported on $\Gamma_a$ can also approximated with measures supported on nearby periodic attractors. It follows that $a \mapsto \mu_a$ is not weak*continuous on any full-measure subset of $(0,2]$. Some of these results extend to unimodal families with critical point of higher order, and even to not-too-flat flat topped families.

## 1. DEFINITIONS, STATEMENTS AND RELATED RESULTS

1.1. **Introduction.** An invariant measure $\mu$ for an interval map $f$ is called a *natural* (*physical*, *Sinai-Ruelle-Bowen*) measure, if $\mu$ describes the asymptotic distribution of $\{f^n(x)\}_{n=0}^{\infty}$ for all $x$ in a set of positive Lebesgue measure, i.e. if

$$\mu \stackrel{\text{weak*}}{=} \lim_{n \to \infty} \frac{1}{n} \sum_{k=0}^{n-1} \delta_{f^k(x)}$$

holds for all $x$ in a positive measure set. Those $x$ for which this do hold, are said to be *generic* for $\mu$.

We consider the quadratic family $f_a(x) = 1 - ax^2$, $0 < a \leq 2$, of dynamical systems on $I = [-1, 1]$. If $f_a$ admits an *acip* (absolutely continuous invariant probability measure) $\mu_a$, then $\mu_a$ is a natural measure for $f_a$, describing the asymptotics of almost all orbits (Theorem V.1.5, [MS93]). If $f_a$ has a periodic attractor $\{x_i\}_{i=1}^p$, then

$$\mu_a := \frac{1}{p} \sum_{i=1}^{p} \delta_{x_i}$$

is again a natural measure for $f_a$ which describes the asymptotics for almost all points ([BL91]).

For a long time it was conjectured that for almost all $a$, $f_a$ either has a periodic attractor or admits an acip. This fundamental result was recently proved by M. Lyubich ([Lyu97]). It follows that the mapping

$$\Psi : a \mapsto \mu_a =: \text{the natural measure of } f_a$$

is well defined for almost all $a$ in $(0, 2]$. Of course, $\Psi$ is continuous on the set of hyperbolic attracting maps. In what follows we discuss the structure in parameter space and the parameter dependence of $\mu_a$ near certain maps admitting an acip, more precisely the maps $f_a$ considered by Benedicks and Carleson.

**Theorem** (Benedicks and Carleson, [BC85], [BC91])**.** *There exists a set $\mathcal{A} \subset (0, 2]$ of positive Lebesgue measure, with $2$ a Lebesgue density point of $\mathcal{A}$, such that if $a \in \mathcal{A}$, then the Lyapunov exponent of $f_a$ at the critical value is positive and also $f_a$ admits an acip $\mu_a$ with a density that belongs to $L^p$ for any $1 \leq p < 2$. Furthermore, for almost all $a \in \mathcal{A}$, the critical point $c = 0$ is generic for $\mu_a$.*

Jakobson, in [Jak81], was the first to prove the existence of a positive measure set of parameters admitting acips. Various versions and generalizations can be found in [Ryc88], [BY92], [TTY92], [Tsu93b], [Tsu93a], [MS93], [Thu97] and [Luz98].

1.2. **Theorems and corollaries.**

**Theorem A.** *For each $a \in \mathcal{A}$, there exists a sequence $\{a_n\}_{n=1}^{\infty}$ such that $f_{a_n}$ has a super-stable periodic attractor of length $r_n$, such that*
  (i) $a_n \to a$, $n \to \infty$;
  (ii) $r_n \uparrow \infty$, $n \to \infty$;
  (iii) *if $x = 0$ is generic for $\mu_a$, then $\mu_{a_n} \stackrel{\text{weak}*}{\longrightarrow} \mu_a$, $n \to \infty$.*

In [Thu96] a version of Theorem A is proven, assuming only sub-exponential growth of the derivative along the orbit of the critical value. In [Ure95] and [Ure96], Ures proves corresponding statements for the Hénon family.

**Theorem B.** *Let $\Gamma = \Gamma_a$ be a hyperbolic set for $f_a$, $a \in \mathcal{A}$, let $z = z(a)$ be any point in $\Gamma$ and let $\Gamma_b$ and $z(b)$ be the continuations of $\Gamma$ and $z$. Then*
$$a \in \mathrm{cl}\{b \mid f_b^N(0) = z(b) \text{ for some } N = N(b) \}$$

*Remark* 1. An obvious choice is $\Gamma_a = \{z(a)\}$, where $z(a)$ is the interior unstable fixed point. This implies that (this type of) Misiurewicz points are dense in $\mathcal{A}$. By taking $\mathcal{A}$ sufficiently close to 2, the theorem holds for any hyperbolic set of $f_2$.

**Theorem C.** *Let $\{x_1, x_2, \ldots, x_p\}$ be a hyperbolic periodic repeller for $f_a$, $a \in \mathcal{A}$. Then there exists a sequence of parameters $\{a_n\}_{n=1}^{\infty}$ converging to $a$ such that $f_{a_n}$ has a super-attractor and such that*
$$\mu_{a_n} \xrightarrow{\text{weak}^*} \mu_a^{\text{sing}} =: \frac{1}{p} \sum_{i=1}^{p} \delta_{x_i}.$$

We also state a theorem that holds for any post-critically finite Misiurewicz map in the quadratic family.

**Theorem D.** *Suppose $f_a$ is a quadratic map whose critical orbit is pre-periodic to an unstable periodic orbit $\{x_1, x_2, \ldots, x_p\}$. Then there is a sequence of parameters $\{a_n\}_{n=1}^{\infty}$ accumulating an $a$ such that $f_{a_n}$ has a super-stable period attractor and*
$$\mu_{a_n} \xrightarrow{\text{weak}^*} \mu_a^{\text{sing}} =: \frac{1}{p} \sum_{i=1}^{p} \delta_{x_i}.$$

*Remark* 2. Since $\Psi : a \mapsto \mu_a$ is continuous when restricted to a periodic window, the statements of theorems A, C and D remains true for any sequence of parameters running thru the periodic windows $J_n \ni a_n$.

From Theorem C and Remark 2, we obtain

**Corollary 1.** *$\Psi : a \mapsto \mu_a$ is not continuous at any point in $\mathcal{A}$, and $\Psi$ is not continuous on any full-measure subset of $(0, 2]$.*

Since Misiurewicz maps admits an acip under some mild conditions, fulfilled in the quadratic family, [Mis81] and [BM89], Theorem D gives the following:

**Corollary 2.** *$\Psi : a \mapsto \mu_a$ is discontinuous at every post-critically finite Misiurewicz map in the quadratic family.*

Most of this can be generalized to generic unimodal families with critical point of higher order, and even to those flat-topped families that are considered in [Thu97]. What may happen is that the last statement of Theorem A becomes empty; it is not known if the critical point is generic for a positive set of acip-maps in these cases.

1.3. **Related results.** The measures $\mu_a$, $a \in \mathcal{A}$, are known to be stable under random perturbations of the iterations of $f_a$ ([BY92], [BV96]). Noise is thus in this way a savior, when modeling/experimenting with unimodal maps in the chaotic regime.

In [Ryc88], Rychlik gave a new proof of Jakobsons Theorem. A positive measure set $\mathcal{A}'$ of parameters is constructed such that for $a \in \mathcal{A}'$, $f_a$ admits an acip $\mu_a$ with a density $\nu_a \in L^p$ for $1 \leq p < 2$. Then in [RS92], Rychlik and Sorret among other things proved that $a \mapsto \nu_a$, defined on $\mathcal{A}'$, is continuous in $L^p$, $1 \leq p < 2$, with a Hölder estimate at Misiurewicz points. The set $\mathcal{A}'$ also has full density at $a = 2$, and so it has a fat intersection with the set $\mathcal{A}$ considered in this paper.

As mentioned, Tsujii has a proof of the Benedicks-Carleson-Jakobson Theorem (leading to an a priori different positive measure set $\mathcal{A}''$ of Collet-Eckmann maps), and some generalizations there of, [Tsu93b] and [Tsu93a]. In [Tsu96] he also discusses weak*-continuity properties of the invariant measures, and among other things proves continuity of $\Psi_{|\mathcal{A}''}$ at Misiurewicz points. Also in this case, 2 is a Lebesgue density point of the good set $\mathcal{A}''$. He also constructs a set $F$ of hyperbolic attracting maps accumulating on $a = 2$, who's natural measures converge to a point mass at the unstable fixed point at $x = -1$ as $a$ tends to 2, and shows that $|F \cap [2 - \epsilon, 2]| \gtrsim \epsilon^2$.

We also remark that the special role of the parameter $a = 2$ is not that special, similar statements holds near any post-critically finite Misiurewicz parameter.

The following is also relevant to this discussion, even though it is a measure 0 phenomena: In [HK90] one constructs

— an uncountable set of parameters, accumulating on $a = 2$, such that the corresponding maps do not have any natural measure at all;
— an uncountable set of parameters, accumulating on $a = 2$, such that the corresponding maps have natural measures $\mu_a = \delta_{z(a)}$.

The first examples of maps with no natural measures was given in [Joh87].

## 2. Properties of the set $\mathcal{A}$

In sections 2.1 – 2.3 we recall some basic facts and definitions from the construction of the set $\mathcal{A}$ that will be used in the proofs. We just state the results that we need in the sequel, and statements do not appear in their logical order. For proofs and more details we refer to [BC85], [BC91] and the expositions in [MS93] and [Luz98]. Similar statements with similar proofs can also be found in [Thu97], dealing with a Benedicks-Carleson theorem for certain flat-top families.

### 2.1. A partition on the interval.
A small neighborhood $I^* = (-\delta, \delta)$ of the critical point $x = 0$ is chosen. $I^*$ is partioned into subintervals

$$I^* \setminus \{0\} = \bigcup_{|\mu| \geq -\log \delta} I_\mu = \bigcup_{\substack{|\mu| \geq -\log \delta \\ 1 \leq \nu \leq \mu^2}} I_{\mu\nu},$$

where $I_\mu = \left[e^{-(\mu+1)}, e^{-\mu}\right)$ for $\mu > 0$, $I_\mu = -I_{-\mu}$ for $\mu < 0$, and

$$I_\mu = \bigcup_{1 \leq \nu \leq \mu^2} I_{\mu\nu}$$

is a subdivision of $I_\mu$ into $\mu^2$ intervals of equal length. $I_{\mu\nu}^+$ denotes the union of $I_{\mu\nu}$ and its two nearest neighbors. We also define $\hat{I}_\mu = \left(-e^{-\mu}, e^{-\mu}\right)$.

### 2.2. Partitions in parameter-space and the mappings $\xi_n$.
The set $\mathcal{A}$ is given as

$$\mathcal{A} =: \bigcap_{n \geq 0} A_n,$$

where $A_n$ is a decreasing sequence of sets in a small one-sided neighborhood of $a = 2$. $A_n$ is constructed by deleting from $A_{n-1}$ according to two principles, which together guarantees that for $a \in A_n$,

(2.1) $$\left|Df_a^j(f_a(0))\right| \geq e^{\lambda j}, \quad \forall j \leq n,$$

for some $\lambda > 0$ independent of $n$ and $a$. One of the exclusion principles simply requires that

(2.2) $$|f_a^n(0)| \geq e^{-\alpha n}.$$

for some suitable, small $\alpha > 0$.

The mappings $\xi_n$ from parameter space to dynamical space are defined via

$$\xi_n(a) := f_a^n(0).$$

These mappings are expanding as long as $f_a^n$ is, in the following sense:

**Lemma 1.** *There is a constant $C$ such that for all $a$ sufficiently close to 2 the following holds: If $|D_x f_a^j(1)| \geq e^{\lambda j}$ for all $j \leq k$ then*
$$\frac{1}{C} \leq \left| \frac{D_a \xi_{k+1}(a)}{D_x f_a^k(1)} \right| \leq C.$$

On each $A_n$ there is a partition $\mathcal{P}_n$ into intervals; each $a \in \mathcal{A}$ is given as $a = \bigcap_{n \geq 0} \omega_n(a)$, where $\omega_n(a)$ is the element of $\mathcal{P}_n$ containing $a$. For each $a$ there is a sequence of times $n_k(a)$, *the essential free return times*, with the following properties:

- $n_k \leq n < n_{k+1} \implies \omega_n = \omega_{n_k}$;
- $n \leq n_k, \xi_n(\omega_{n_k}) \cap I^* \neq \emptyset \implies \xi_n(\omega_{n_k}) \subset$ some $I_{\mu\nu}^+$;
- $I_{\mu\nu} \subset \xi_{n_k}(\omega_{n_k}) \subset I_{\mu\nu}^+$, for some $I_{\mu\nu}$.

A free return is followed by a so called *bound period*, when $\xi_{n_k+j}(a) = f_a^{n_k+j}(0)$ shadows an initial segment $f_a^j(0)$ of the critical orbit closely. More precisely: If $I_{\mu\nu} \subset \xi_{n_k}(\omega_{n_k}) \subset I_{\mu\nu}^+$, then $\xi_{n_k+j}(\omega)$ is in a bound period as long as

$$(2.3) \qquad \left| \bigcup_{a \in \omega} f_a^j \left( -e^{-\mu}, e^{-\mu} \right) \right| \leq e^{-2\alpha j}.$$

We use $p = p(\mu, \omega)$ to denote the length of a bound period.

As long as orbiting outside $I^*$, $\xi_n(\omega)$ grows exponentially (Lemma 3). The small derivative picked up at a return to $I^*$ is compensated for during the bound period by an inductive argument, the net effect is in fact a weaker exponential growth. In particular bound periods are always of finite length, (Lemma 2). Let $p_k$ temporally denote the length of the bound period following a return at time $n_k$. By definition, $n_{k+1}(a)$ is the smallest integer $j \geq n_k(a) + p_k(a)$ such that $\xi_j(\omega_{n_k}(a)) \supset$ some $I_{\mu\nu} \subset I^*$. A return to $I^*$ at some time $j$, $n_k + p_k \leq j < n_{k+1}$, when no $I_{\mu\nu}$ is covered, is called an *inessential free return*. Such returns are also followed by a bound period, after which comes a free period terminating in a new free return. It can be shown that an essential free return always occur after finitely many steps. The "dynamics" of $\{\xi_n\}$ is described more precisely in the following lemmas.

**Lemma 2.** *Suppose inequality (2.1) and condition (2.2) hold for all $j \leq n$ and all $a \in \omega$, and suppose that $\xi_n(\omega) \subset I_{\mu\nu}^+$ where $|\mu| \geq -\log \delta$ and $\delta$ is sufficiently small. Then there exists positive constants $C_0$ and $C$, independent of $\delta$, such that*

(i) *for all $a \in \omega$, all $y \in f_a(\hat{I}_\mu)$ and all $j \leq p$,*
$$\frac{1}{C_0} \leq \frac{|D_x f_a^j(y)|}{|D_x f_a^j(1)|} \leq C_0;$$

(ii) $C\left|\mu\right| \leq p(\mu,\omega) \leq 3\left|\mu\right|/\lambda \leq 3\alpha n/\lambda < n/100$;
(iii) $\left|(f_a^p)'(x)\right| \geq Ce^{\lambda p/4}, \quad \forall x \in I_\mu$;
(iv) $\left|\xi_{n+p(\mu,\omega)}(\omega)\right| \geq Ce^{\lambda p/4}\left|\xi_n(\omega)\right|$;

**Lemma 3.** *There is a $\lambda_0 > 0$ such that for any small $\delta > 0$ and $a_0$ sufficiently close to 2 the following holds: Suppose that $\omega \subset [a_0, 2]$ is such that $\xi_{\hat{n}}(\omega) \subset I_{\mu\nu}$ and $\xi_n(\omega)$ are two consecutive free returns with return times $\hat{n}$ and $n$, $\hat{n} < n$. Also assume that $|Df_a^j(1)| \geq e^{\lambda j}$ for all $j < n$ and all $a \in \omega$. Then there is a constant $C$, independent of $\delta$, such that the following holds:*

(i) $\quad |\xi_{n-k}(\omega)| \leq Ce^{-\lambda_0 k}|\xi_n(\omega)|, \quad \forall 1 \leq k \leq n - \hat{n} - p(\mu, \omega)$;
(ii) $\quad |\xi_n(\omega)| \geq 2|\xi_{\hat{n}}(\omega)|.$

*Furthermore there is a positive integer $N_0(\delta)$ such that for any $\omega$ close to 2,*

(iii) $\quad \xi_{k+j}(\omega) \cap I^* = \emptyset, \quad j = 0, 1, \ldots, N_0 \implies$
$$|\xi_{k+N_0}(\omega)| \geq e^{\frac{\lambda_0}{2}N_0}|\xi_k(\omega)|.$$

Let HD-dist$(J, K)$ denote the Hausdorff distance between the sets $J$ and $K$.

**Lemma 4.** *Suppose inequality (2.1) and condition (2.2) hold for all $j \leq n$ and all $a \in \omega$, and suppose that $\xi_n(\omega) \subset I_\mu$ where $|\mu| \geq -\log \delta$. If $\omega$ is sufficiently close to 2, then for each $a, b \in \omega$ we have*

(i) $\quad$ HD-dist$\left(f_a^j\left(\hat{I}_\mu\right), f_b^j\left(\hat{I}_\mu\right)\right) < \dfrac{1}{1000}\left|f_a^j\left(\hat{I}_\mu\right)\right|$

(ii) $\quad$ HD-dist$\left(\xi_{n+j+1}(\omega), f_a^j\left(\xi_{n+1}(\omega)\right)\right) < \dfrac{1}{1000}\left|f_a^j\left(\xi_{n+1}(\omega)\right)\right|$

*for all $j \leq p(\mu, \omega)$.*

2.3. **Escape-times.** An important role is played by the *escape-times* of $a$. They are an infinite subsequence of $\{n_k(a)\}$, defined via the condition

$n_k$ is an escape-time

(2.4) $\qquad\qquad\qquad \Leftrightarrow$

$\xi_{n_k}(\omega_{n_k-1})$ intersects $(-\delta^2, \delta^2)$ and $\left|\xi_{n_k}(\omega_{n_k-1})\right| \geq \delta$.

It will be important that each $a \in \mathcal{A}$ experiences infinitely many escape-times. This is a consequence of the second parameter selection principle, which roughly speaking discards parameters that on the average has to wait too long for their escape-situations.

## 3. Escaping

We prepare the proofs of theorems A and B by showing that any escaping parameter interval $\omega_n(a)$ contains a super-stable parameter as well as parameters for which the critical orbit lands on any prescribed point in any hyperbolic set of $f_a$.

**Lemma 5.** *Pick an $a \in \mathcal{A}$ and a hyperbolic set $\Gamma_a$ for $f_a$ and a point $z(a) \in \Gamma_a$. Let $n_k$ by an escape time for $a$. Then there are two parameters $a^*, \hat{a} \in \omega_{n_{k-1}}(a)$ and two integers $r^*$ and $\hat{r}$, $0 < r^* \leq \hat{r} \lesssim -\log \delta$ such that*

(i) *$f_{a^*}$ has a super-stable attractor of period $n_k(a) + r^*$;*
(ii) *$f_{\hat{a}}^j(0) \neq z(\hat{a})$ for $j < n_k + \hat{r}$ and $f_{\hat{a}}^{n_k+\hat{r}}(0) = z(\hat{a})$.*

*Proof.* Let $\omega_{n_{k-1}}(a) = (b, c)$. We may assume that $\xi_{n_k}(b) = \delta^2$ and $\xi_{n_k}(c) = \delta$. The idea is that distance between $\xi_{n_k+j}(b)$ and $-1$ will be $\leq 4^j \delta^4$ for $j > 1$, while the distance between $-1$ and $\xi_{n_k+j}(c)$ will be $\geq 3^j \delta^2$ as long as $\xi_{n_k+j}(c) < -3/4$. From this it follows that for some $j_0 \lesssim -\log \delta$, $\xi_{n_k+j_0}(\omega_{n_k})$ will grow to length $\sim 1/4$, with its left end $\xi_{n_k+j_0}(b)$ still within a distance $o(\delta)$ from $-1$.

Thus $\xi_{n_k+r^*}$ maps $(b, c)$ across $x = 0$ for $r^* = j_0 + 1$ or $r^* = j_0 + 2$, and since $0 \notin \xi_j(b, c)$ for $j < n_k + j^*$ by (2.2), the required $a^*$ has been found.

Since $\Gamma_a$ is hyperbolic, $z(a)$ moves continuously with $a$, and $(z(b), z(c))$ will be a very small interval. It follows that $\xi_{n_k+\hat{r}}(b, c) \supset (z(b), z(c))$ for $\hat{r} = r^*$ or $\hat{r} = r^* + 1$. Thus $\xi_{n_k+\hat{r}}(a) - z(a)$ changes sign on $(b, c)$, and so $\xi_{n_k+\hat{r}}(\hat{a}) = z(\hat{a})$ for some $\hat{a} \in (b, c)$. $\square$

## 4. Proof of Theorem A and Theorem B

Each $a \in \mathcal{A}$ is given as

$$a = \cap_{n=1}^{\infty} \omega_n, \quad \omega_n \in \mathcal{P}_n,$$

with infinitely many escape situations where Lemma 5 can be applied. For any $b \in \omega_n$, inequality (2.1) holds. Thus we may use Lemma 1 to conclude that $|\omega_n| \lesssim e^{-n\lambda}$. From this the first two statements of Theorem A and Theorem B follows.

We now prove the last part of Theorem A. Let $a \in \mathcal{A}$ be such that $x = 0$ is generic for $\mu_a$, and let $\{a_n\}$ be the sequence of super-stable parameters converging to $a$ constructed above.

We have to show that

$$\lim_{n \to \infty} \int \varphi \, d\mu_{a_n} = \int \varphi \, d\mu_a$$

for any continuous function $\varphi$. It is enough to consider Lipschitz continuous test functions. Let $\kappa$ be the Lipschitz' constant of $\varphi$, and let $r_k$ denote the length of super-stable attractor of $f_{a_k}$. We have that for any $a_k$,

$$\left| \int \varphi \, d\mu_a - \int \varphi \, d\mu_{a_k} \right|$$
$$\leq \left| \int \varphi \, d\mu_a - \frac{1}{r_k} \sum_{j=0}^{r_k-1} \varphi(\xi_j(a)) \right|$$
$$+ \left| \frac{1}{r_k} \sum_{j=0}^{r_k-1} (\varphi(\xi_j(a)) - \varphi(\xi_j(a_k))) \right|$$
$$+ \left| \frac{1}{r_k} \sum_{j=0}^{r_k-1} \varphi(\xi_j(a_k)) - \int \varphi \, d\mu_{a_k} \right|,$$

where the last term $= 0$ by definition, and the first one is $< \epsilon$ for all $k > N_0(\epsilon)$ since 0 is generic for $\mu_a$ and since $r_k \uparrow \infty$ when $k \to \infty$. The second term finally is

$$\leq \frac{\kappa}{r_k} \sum_{j=0}^{r_k-1} |\xi_j(a) - \xi_j(a_k)|,$$

so it suffices to prove $S := \sum_{j=0}^{r_k-1} |\xi_j(a) - \xi_j(a_k)|$ is bounded by some constant independent of $a_k$. Let

$$\Lambda_j = |\xi_j(a) - \xi_j(a_k)|,$$

and let $n_k$ be the escape-time preceding the creation of the super-stable orbit at time $r_k$; then $r_k - n_k \lesssim \log 1/\delta$ (c.f. Lemma 5). Remember that $a_k \in \omega_{n_k}(a) \in \mathcal{P}_{n_k}$, and that for the "orbit" of $\omega_{n_k}$ under the family $\{\xi_i\}_{i=1}^{n_k}$ certain free return times $\{t_i\}_{i=1}^{T}$ are defined, $t_T = n_k$, and that each free return $t_i$ is followed by a bound period of finite length, which we denote $p_i$. For the sake of notation, we define $t_0 = p_0 = 0$. We split the sum $S$ into sub-sums:

$$S := \sum_{i=0}^{T-1} \left( S_i^{\text{bp}} + S_i^{\text{fp}} \right) + S^{\text{tail}},$$

where

$$S_i^{\text{bp}} = \sum_{l=t_i}^{t_i+p_i-1} \Lambda_l, \qquad S_i^{\text{fp}} = \sum_{l=t_i+p_i}^{t_{i+1}-1} \Lambda_l$$

and

$$S^{\text{tail}} = \sum_{l=t_T}^{r_k} \Lambda_l.$$

With this notation, $S_0^{\text{fp}}$ is the contribution up till the first free return and $S_0^{\text{bp}}$ is an empty sum and equals 0. First we observe that $S^{\text{tail}}$ has no more than $\sim -\log \delta$ terms, and is therefore $\leq C(\delta)$. We now estimate $\sum_{i=0}^{T-1} S_i^{\text{fp}}$, using (i) and (ii) of Lemma 3:

$$\sum_{i=0}^{T-1} S_i^{\text{fp}} = \sum_{i=0}^{T-1} \sum_{l=t_i+p_i}^{t_{i+1}-1} \Lambda_l \leq \sum_{i=0}^{T-1} \sum_{l=t_i+p_i}^{t_{i+1}-1} Ce^{-\lambda(t_{i+1}-l)} \Lambda_{t_{i+1}}$$

$$\leq C_1 \sum_{i=1}^{T} \Lambda_{t_i} \leq C_1 \sum_{i=1}^{T} 2^{i-T} \Lambda_{t_T} \leq C_2.$$

We now turn to $\sum_{i=0}^{T-1} S_i^{\text{bp}}$. First we estimate the individual terms in $S_i^{\text{bp}}$. Now $\Lambda_l = |\xi_l(a) - \xi_l(a_k)|$, and $(a; a_k) \subset \omega_{n_k} \in \mathcal{P}_{n_k}$. Since $t_i$ is a free return, $t_i < t_T = n_k$, it follows that $\xi_{t_i}((a; a_k)) \subset$ some $I_{\mu_i \nu_i}^+ \subset I^*$. Using Lemma 2, Lemma 4 and the binding condition (2.3), we see that for $1 \leq j < p_i$ and any $b \in (a; a_k)$,

$$\Lambda_{t_i+j} \lesssim \left|f_b^j\left(\xi_{t_i}((a; a_k))\right)\right| = \frac{\left|f_b^j\left(\xi_{t_i}((a; a_k))\right)\right|}{\left|f_b^j\left(\hat{I}_{\mu_i}\right)\right|} \left|f_b^j\left(\hat{I}_{\mu_i}\right)\right|$$

$$\lesssim \frac{|f_b(\xi_{t_i}((a; a_k)))|}{\left|f_b\left(\hat{I}_{\mu_i}\right)\right|} e^{-2\alpha j} \sim \frac{|f_b(\xi_{t_i}((a; a_k)))|}{|f_b(I_{\mu_i})|} e^{-2\alpha j}$$

$$\sim \frac{|(\xi_{t_i}((a; a_k)))|}{|I_{\mu_i}|} e^{-2\alpha j}$$

In the last two steps we also used the facts that $\left|f_b\left(\hat{I}_\mu\right)\right| \sim |f_b(I_\mu)|$, and that the distortion of $f_b$ restricted to $I_\mu$ has a bound independent of $\mu$.

Obviously $\Lambda_{t_j} < |\xi_{t_i}((a; a_k))|/|I_{\mu_i}|$, so we obtain

$$S_i^{\text{bp}} \lesssim \sum_{j=0}^{\infty} \frac{|(\xi_{t_i}((a; a_k)))|}{|I_{\mu_i}|} e^{-2\alpha j} \lesssim \frac{|(\xi_{t_i}((a; a_k)))|}{|I_{\mu_i}|} \leq 3 \frac{|I_{\mu_i \nu_i}|}{|I_{\mu_i}|} \lesssim \frac{1}{\mu_i^2}.$$

Finally we estimate $\sum S_i^{\text{bp}}$:

$$\sum_{i=0}^{T-1} S_i^{\text{bp}} \lesssim \sum_{\mu_i} \sum_{\substack{\text{returns} \\ \text{to } I_{\mu_i}}} \mu_i^{-2} \lesssim \sum_{\mu_i} \sum_{\substack{\text{last return} \\ \text{to } I_{\mu_i}}} \mu_i^{-2}$$

$$\lesssim \sum_{|\mu| \geq -\log \delta} \mu^{-2} < C.$$

The second "$\lesssim$" holds because of (ii) of Lemma 3. This finishes the proof of Theorem A.

## 5. Proof of Theorem D

Suppose $\Gamma_a = \{x_1, x_2, \ldots, x_p\}$ is a hyperbolic repelling periodic orbit for $f_a$, absorbing the critical orbit. Let $N$ be minimal such that $f_a^N(0) \in \Gamma_a$, and assume that $f_a^N(0) = x_1$. For $\gamma > 0$ we define $J_\gamma = \bigcup_{i=1}^p [x_i - \gamma, x_i + \gamma]$.

Since $f_a$ is Misiurewicz map with finite critical orbit, we can carry out the Benedicks-Carleson construction in a neighborhood of $a$ and find a sequence $\{\omega_n\}_{n=N}^\infty$ of parameter intervals such that

- $\omega_{n+1} \subset \omega_n$ and $a = \bigcap_{n=N}^\infty \omega_n$;
- $\xi_k(\omega_n) \subset J_\gamma$ for $k = N, \ldots, n$;
- for all $n \geq N$ there is an $i_n$ such that $\xi_n(\omega_n) = [x_{i_n} - \gamma, x_{i_n} + \gamma]$;
- there is a natural number $\hat{N} = \hat{N}(\gamma)$ such that for each $n \geq N$ there is an $m_n \leq \hat{N}$ such that $\xi_{n+m_n}(\omega_n) \ni 0$.

This follows from continuity, uniform expansion away from the critical point and the fact that parameter- and space-derivatives are comparable. It follows that for all $n \geq N$, there is an $a_n \in \omega_n$ such that $f_{a_n}$ has super-attractor of length $n + m_n$ such that

$$\#\{j \leq n + m_n \mid f_{a_n}^j(0) \notin J_\gamma\} = N - 1 + m_n \leq N + \hat{N}(\gamma).$$

Now pick a $\varphi \in C^0(I)$ and an $\epsilon > 0$. Since $\varphi$ is uniformly continuous, we can fix a $\gamma$ such that $|x - y| < \gamma$ implies $|\varphi(x) - \varphi(y)| < \epsilon/2$. If $\mu_{a_n}$ is the natural measure for $f_{a_n}$, where $a_n$ is as above, we have that

$$\int \varphi \, d\mu_{a_n} = \frac{1}{n + m_n} \sum_{i=1}^{n+m_n} \varphi\left(f_{a_n}^i(0)\right)$$

$$\leq \frac{1}{n + m_n} \sum_{i=N}^{n} \varphi\left(f_{a_n}^i(0)\right) + \frac{N - 1 + m_n}{n + m_n} \sup_I \varphi$$

$$\leq \frac{1}{n + m_n} \sum_{i=N}^{n} \varphi\left(f_{a_n}^i(0)\right) + \frac{N + \hat{N}}{n + \hat{N}} \sup_I \varphi.$$

It is clear that the last term vanishes when $n \to \infty$. We have to show that $\frac{1}{n+m_n} \sum_{i=N}^{n} \varphi\left(f_{a_n}^i(0)\right)$ tends to $\int \varphi \, d\mu_a^{\text{sing}} = \frac{1}{p} \sum_{i=1}^{p} \varphi(x_i)$, when $n$ tends to infinity. Without loss of generality we may assume that $n - N + 1 = R_n p$ for some integer $R_n$. Let $y_i = f_{a_n}^{N+i-1}(0)$ for $1 \leq i \leq n + 1 - N$, and for $i > p$ define $x_i = f_a^{i-1}(x_1)$. Then

$$\frac{1}{n+m_n} \sum_{i=N}^{n} \varphi\left(f_{a_n}^i(0)\right)$$

$$= \frac{1}{n+m_n} \sum_{i=1}^{n+1-N} \varphi(y_i)$$

$$\leq \frac{pR_n}{n+m_n} \frac{1}{pR_n} \sum_{i=1}^{n+1-N} (\varphi(x_i) + \epsilon/2)$$

$$= \frac{n-N+1}{n+m_n} \left( \int \varphi \, d\mu_a^{\text{sing}} + \epsilon/2 \right)$$

$$\leq \int \varphi \, d\mu_a^{\text{sing}} + \epsilon,$$

where the last inequality holds for $n$ sufficiently large. In the same way

$$\int \varphi \, d\mu_{a_n} \geq \int \varphi \, d\mu_a^{\text{sing}} - \epsilon$$

for $n$ sufficiently large.

## 6. Proof of Theorem C

Since $\{x_1, \ldots, x_p\}$ is a hyperbolic repeller, it persists for nearby parameter values. Let $\Gamma_a = \{x_1(a), \ldots, x_p(a)\}$ by its continuation. By Theorem B, we know that there is a sequence $\{b_n\}$ converging to $a$, such that $f_{b_n}$ has critical point pre-periodic to $x_1(b_n)$. By Theorem D, for each $b_n$ there is a sequence $\{a_{n,m}\}_{m=1}^{\infty}$ corresponding to maps with super stable attractors, such that $\lim_{m \to \infty} a_{n,m} = b_n$ and $\lim_{m \to \infty} \mu_{a_{n,m}} = \sum_{i=1}^{p} \delta_{x_i(b_n)}$. Since $\Gamma_a$ moves continuously with $a$, we get a sequence $\{a_{n,m(n)}\}_{n=1}^{\infty}$, corresponding to super attractors and converging to $a$ such that

$$\lim_{n \to \infty} \mu_{a_{n,m(n)}} \stackrel{\text{weak}^*}{=} \sum_{i=1}^{p} \delta_{x_i(a)}.$$

## References

[BC85]   M. Benedicks and L. Carleson, *On iterations of* $1 - ax^2$, Ann. of Math. **122** (1985), 1–25.

[BC91]  M. Benedicks and L. Carleson, *The dynamics of the Hénon map*, Ann. of Math. **133** (1991), 73–169.
[BL91]  A.M. Blokh and M.Yu. Lyubich, *Measurable dynamics of S-unimodal maps of the interval*, Ann. Sc. E.N.S. 4e série **24** (1991), 545–573.
[BM89]  M. Benedicks and M. Misiurewicz, *Absolutely continuous invariant measures for maps with flat tops*, Inst. Hautes Études Sci. Publ. Math. (1989), no. 69, 203–213.
[BV96]  V. Baladi and M. Viana, *Strong stochastic stability and rate of mixing for unimodal maps*, Ann. Sci. cole Norm. Sup. (4) **29** (1996), no. 4, 483–517.
[BY92]  M. Benedicks and L. S. Young, *Absolutely continuous invariant measures and random perturbations for certain one-dimensional maps*, Ergod. Th. & Dynam. Sys. **12** (1992), 13–37.
[HK90]  F. Hofbauer and G. Keller, *Quadratic maps without asymptotic measure*, Comm. Math. Phys. **127** (1990), no. 2, 319–337.
[Jak81]  M. Jakobson, *Absolutely continuous invariant measures for one-parameter families of one-dimensional maps*, Comm. Math. Phys. **81** (1981), 39–88.
[Joh87]  S.D. Johnson, *Singular measures without restrictive intervals*, Comm. Math. Phys. **110** (1987), no. 2, 185–190.
[Luz98]  S. Luzzatto, *A proof of Jakobson's theorem*, In preparation, 1998.
[Lyu97]  M. Lyubich, *Almost every real quadratic map is either regular or stochastic*, SUNY preprint IMS 97–8, SUNY Stony-Brook, New York, 1997.
[Mis81]  M. Misiurewicz, *Absolutely continuous measures for certain maps of an interval*, Inst. Hautes Études Sci. Publ. Math. (1981), no. 53, 17–51.
[MS93]  W. de Melo and S. van Strien, *One-dimensional dynamics*, Springer-Verlag, Berlin-Heidelberg, 1993.
[RS92]  M. Rychlik and E. Sorets, *Regularity and other properties of absolutely continuous invariant measures for the quadratic family*, Comm. Math. Phys. **150** (1992), no. 2, 217–236.
[Ryc88]  M. Rychlik, *Another proof of Jakobsons theorem and related results*, Ergod. Th. & Dynam. Sys. **8** (1988), 93–109.
[Thu96]  H. Thunberg, *Absolutely continuous invariant measures and super-stable periodic orbits: weak\*-convergence of natural measures*, Some Problems in Unimodal Dynamics (Stockholm), KTH, Stockholm, 1996, Ph.D. thesis, pp. I: 1–19.
[Thu97]  H. Thunberg, *Positive exponent in families with flat critical point*, To appear in Ergod. Th. and Dynam. Sys., 1997.
[Tsu93a]  M. Tsujii, *Positive Lyapunov exponent in families of one-dimensional dynamical systems*, Inven. math. **111** (1993), no. 1, 113–137.
[Tsu93b]  M. Tsujii, *A proof of Benedicks-Carleson-Yakobson theorem*, Tokyo J. Math. **16** (1993), no. 2, 295–310.
[Tsu96]  M. Tsujii, *On continuity of Bowen-Ruelle-Sinai measures in families of one-dimensional maps*, Comm. Math. Phys. **117** (1996), no. 1, 1–11.
[TTY92]  P. Thieullen, C. Tresser, and L.S. Young, *Exposant de Lyapunov positif dans des familles à un paramètre d'applications unimodales*, C. R. Acad. Sci. Paris, Série I **t315** (1992), 69–72.
[Ure95]  R. Ures, *On the approximation of Hénon-like attractors by homoclinic tangencies*, Ergodic Theory Dynam. Systems **15** (1995), no. 6, 1223–1229.

[Ure96]  R. Ures, *Hénon attractors: SBR measures and Dirac measures for sinks*, International Conference on Dynamical Systems (Montevideo, 1995) (Harlow), Pitman Res. Notes Math. Ser., vol. 362, Longman, Harlow, 1996, pp. 214–219.

Department of Mathematics, KTH, S-100 44 STOCKHOLM, Sweden
*E-mail address*: `hasset@math.kth.se`